\documentclass[12pt]{article}
\usepackage{amsmath}
\usepackage{amsfonts}

\usepackage{amsthm} 
\usepackage{amscd}
\theoremstyle{plain} 
\newtheorem{theorem}{Theorem}[section]

\newtheorem{lemma}[theorem]{Lemma}

\numberwithin{equation}{section}

\newcommand{\al}{\alpha}

\newcommand{\br}{\mathbf {R}}

\newcommand{\x}{\mathcal{X}}
\newcommand{\om}{\omega}
\newcommand{\ot}{\otimes}
\newcommand{\w}{\wedge}

\newcommand{\p} {\partial}
\newcommand{\g} {\mathfrak {g}}
\newcommand{\s}{\mathfrak {sp}}

\begin{document}
\title{A Characteristic Map for Symplectic Manifolds}
\author{Jerry M. Lodder}
\date{}
\maketitle

\noindent
{\em Mathematical Sciences, Dept. 3MB  \\
Box 30001 \\ 
New Mexico State University \\
Las Cruces NM, 88003, U.S.A. }

\noindent
e-mail:  {\em jlodder@nmsu.edu}

\bigskip
\noindent
{\bf Abstract.}  We construct a local characteristic map to a
symplectic manifold $M$ via certain cohomology groups of Hamiltonian
vector fields.  For each $p \in M$, the Leibniz cohomology of the
Hamiltonian vector fields on $\br^{2n}$ maps to the Leibniz cohomology
of all Hamiltonian vector fields on $M$.  For a particular extension
$\g_n$ of the symplectic Lie algebra, the Leibniz cohomology of $\g_n$
is shown to be an exterior algebra on the canonical symplectic two-form.  The
Leibniz homology of $\g_n$ then maps to the Leibniz homology of
Hamiltonian vector fields on $\br^{2n}$.   

\bigskip
\noindent
{\bf Mathematics Subject Classifications (2000):}  17B56, 53D05, 17A32. 

\bigskip
\noindent
{\bf Key Words:}  Symplectic topology, Leibniz homology, symplectic
invariants.

\section{Introduction}

We construct a local characteristic map to a symplectic manifold $M$
via certain cohomology groups of Hamiltonian vector fields.  Recall
that the group of affine symplectomorphisms, i.e., the affine
symplectic group $ASp_n$, is given by all transformations $\psi :
\br^{2n} \to \br^{2n}$ of the form
$$  \psi(z) = Az + z_0 ,  $$
where $A$ is a $2n \times 2n$ symplectic matrix and $z_0$ a fixed
element of $\br^{2n}$ \cite[p.\ 55]{MS}.  Let $\g_n$ denote the Lie
algebra of $ASp_n$, referred to as the affine symplectic Lie algebra.  
Then $\g_n$ is the largest
finite dimensional Lie subalgebra of the Hamiltonian vector fields on 
$\br^{2n}$, and serves as our point of departure for calculations.
Particular attention is devoted to the Leibniz homology 
of $\g_n$, i.e., $HL_*(\g_n ;\, \br )$, and proven is that
$$  HL_*(\g_n ; \, \br ) \simeq \Lambda^* (\om_n), $$
where $\om_n = \sum_{i=1}^n \frac{\p}{\p x^i} \wedge \frac{\p}{\p y^i}$ 
and $\Lambda^*$ denotes the exterior algebra.  Dually, for cohomology,
$$  HL^* (\g_n ; \, \br) \simeq \Lambda^*( \om_n^*), $$
where $\om_n^* = \sum_{i=1}^n dx^i \wedge dy^i$.  

For $p \in M$, the local characteristic map acquires the form
$$  HL^*( \x_{H}(\br^{2n});\, \br )(p) \to HL^* (\x_{H}(M); \,
C^{\infty}(M)), $$
where $\x_{H}$ denotes the Lie algebra of Hamiltonian vector fields,
and $C^{\infty}(M)$ is the ring of $C^{\infty}$ real-valued functions
on $M$.  
Using previous work of the author \cite{Lodder}, there is a natural
map
$$  H_{dR}^*(M; \, \br) \to HL^*( \x(M); \, C^{\infty}(M)), $$
where $H^*_{dR}$ denotes deRham cohomology.  Composing with
$$ HL^*(\x(M); \, C^{\infty}(M))  \to HL^*(\x_H(M); \, C^{\infty}(M)), $$ 
we have
$$ H_{dR}^* (M; \, \br) \to HL^*(\x_{H}(M); \, C^{\infty}(M)). $$
The inclusion of Lie algebras $\g_n \hookrightarrow \x_H (\br^{2n})$
induces a linear map
$$  HL_*(\g_n ; \, \br) \to HL_*(\x_{H}(\br^{2n}); \, \br)  $$
and $HL^*(\x_{H}(\br^{2n}); \, \br)$ contains a copy of 
$HL^*(\g_n ; \, \br)$ as a direct summand.

The calculational tools for $HL_* (\g_n)$ include the Hochschild-Serre
spectral sequence for Lie-algebra (co)homology, the Pirashvili spectral
sequence for Leibniz homology, and the identification of certain
symplectic invariants of $\g_n$ which appear in the appendix.

\section{The Affine Symplectic Lie Algebra}

As a point of departure, consider a $C^{\infty}$ Hamiltonian function
$H: \br^{2n} \to \br$ with the associated Hamiltonian vector field
$$ X_H = \sum_{i=1}^{n} \frac{\p H}{\p x_i} \frac{\p}{\p y^i} -
\sum_{i=1}^{n} \frac{\p H}{\p y_i}\frac{\p}{\p x^i}, $$
where $\br^{2n}$ is given coordinates
$$ \{ x_1, \ x_2, \ \ldots, \ x_n, \ y_1,\ y_2, \ \ldots, \ y_n \}, $$
and $\frac{\p}{\p x^i}$, $\frac{\p}{\p y^i}$ are the unit vector
fields parallel to the $x_i$ and $y_i$ axes respectively.  The vector
field $X_H$ is then tangent to the level curves (or hyper-surfaces) of
$H$.  Restricting $H$ to a quadratic function in  
$$ \{ x_1, \ x_2, \ \ldots, \ x_n, \ y_1,\ y_2, \ \ldots, \ y_n \}, $$
yields a family of vector fields isomorphic to the real symplectic Lie
algebra $\s_n$.  For $i$, $j$, $k \in \{ 1, \ 2, \ 3, \ \ldots ,\
n\}$, an $\br$-vector space basis, $\mathcal{B}_1$, for $\s_n$ is given
by the families: 
\begin{itemize}
\item[(1)]  $x_k \frac{\p}{\p y^k}$ 
\item[(2)]  $y_k \frac{\p}{\p x^k}$
\item[(3)]  $x_i \frac{\p}{\p y^j} + x_j \frac{\p}{\p y^i}$, $i \neq j$
\item[(4)]  $y_i \frac{\p}{\p x^j} + y_j \frac{\p}{\p x^i}$, $i \neq j$
\item[(5)]  $y_j \frac{\p}{\p y^i} - x_i \frac{\p}{\p x^j}$, ($i = j$ possible).
\end{itemize}
It follows that ${\text{dim}}_{\br} (\s_n ) = 2n^2 + n$.

Let $I_n$ denote the Abelian Lie algebra of Hamiltonian vector fields
arising from the linear (affine) functions $H: \br^{2n} \to \br$.
Then $I_n$ has an $\br$-vector space basis given by
$$ {\mathcal{B}}_2 = \Big\{ \frac{\p}{\p x^1}, \ \frac{\p}{\p x^2}, \ 
\ldots, \ \frac{\p}{\p x^n}, \ \frac{\p}{\p y^1}, \ \frac{\p}{\p y^2},
\ \ldots, \ \frac{\p}{\p y^n} \Big\}.  $$
The affine symplectic Lie algebra, $\g_n$, has an $\br$-vector space
basis ${\mathcal{B}}_1 \cup {\mathcal{B}}_2$ .  There is a short exact
sequence of Lie algebras
$$ \CD
0 @>>> I_n @>i>> \g_n @>{\pi}>> \s_n @>>> 0, 
\endCD $$
where $i$ is the inclusion map and $\pi$ is the projection
$$ \g_n \to ( \g_n / I_n ) \simeq \s_n . $$
In fact, $I_n$ is an Abelian ideal of $\g_n$ with $I_n$ acting on
$\g_n$ via the bracket of vector fields.

\section{The Lie Algebra Homology of $\g_n$}

For any Lie algebra $\g$ over a ring $k$, the Lie algebra homology of
$\g$, written $H^{\text{Lie}}_*(\g ; \, k)$, is the homology of the
chain complex $\Lambda^* (\g )$, namely
$$ \CD 
k @<0<< \g @<[\ , \ ]<< \g^{\wedge 2} @<<< \ldots @<<< \g^{\wedge (n-1)}
@<d<< \g^{\wedge n} @<<< \ldots,  
\endCD $$
where
\begin{align*}
& d(g_1 \wedge g_2 \wedge \, \ldots \, \wedge g_n) = \\
& \sum_{1 \leq i < j \leq n} (-1)^j \, (g_1 \wedge \, \ldots \, \wedge
g_{i-1} \wedge [g_i, \, g_j] \wedge g_{i+1} \wedge \, \ldots \, 
\hat{g}_j \,\ldots \, \wedge g_n). 
\end{align*}
For actual calculations in this paper, $k = \br$.  Additionally, Lie
algebra homology with coefficients in the adjoint representation,
written $H^{\text{Lie}}_*(\g ; \, \g )$, is the homology of the chain
complex $\g \otimes \Lambda^*( \g )$, i.e., 
$$ \g \longleftarrow \g \ot \g \longleftarrow \g \ot \g^{\wedge 2} \longleftarrow
\ldots \longleftarrow \g \ot \g^{\wedge (n-1)} \,
{\overset{d}{\longleftarrow}} \,\g \ot \g^{\wedge n} \longleftarrow \ldots ,$$
where
\begin{align*}
& d(g_1 \ot g_2 \wedge g_3 \, \ldots \, \wedge g_{n+1}) = 
\sum_{i=2}^{n+1} (-1)^i \, ([g_1, \, g_i] \ot g_2 \wedge \, \ldots
  \, \hat{g}_i \, \ldots \, \wedge g_{n+1}) \\
& + \sum_{2 \leq i < j \leq n+1} (-1)^j \, (g_1 \ot g_2 \wedge \, \ldots \, \wedge
g_{i-1} \wedge [g_i, \, g_j] \wedge g_{i+1} \wedge \, \ldots \, 
\hat{g}_j \,\ldots \, \wedge g_{n+1}). 
\end{align*}
The canonical projection $\g \ot \Lambda^*(\g ) \to \Lambda^{*+1}(\g )$
given by $\g \ot \g^{\wedge n} \to \g^{\wedge (n+1)}$ is a map of
chain complexes and induces a $k$-linear map on homology
$$  H^{\text{Lie}}_n(\g ; \, \g) \to H^{\text{Lie}}_{n+1}(\g ; \,k). $$

Given a (right) $\g$-module $M$, the module of invariants $M^{\g}$ is
defined as
$$  M^{\g} = \{ m \in M \ | \ [m, \, g] = 0 \ \ \forall g \in \g \}. $$
Note that $\s_n$ acts on $I_n$ and on the affine symplectic Lie
algebra $\g_n$ via the bracket of vector fields.  The action is
extended to $I^{\wedge k}_n$ by 
$$  [\al_1 \wedge \al_2 \wedge \, \ldots \, \wedge \al_k, \, X] =
\sum_{i=1}^k \al_1 \wedge \al_2 \wedge \, \ldots \, \wedge [\al_i, \, X] \wedge \,
\ldots \, \wedge \al_k $$
for $\al_i \in I_n$, $X \in \s_n$, and similarly for the $\s_n$
action on $\g_n \ot I^{\wedge k}_n$.  The main result of this section
of the following.
\begin{lemma} \label{3.1}
There are natural vector space isomorphisms
\begin{align*}
& H^{\text{Lie}}_*(\g_n ; \, \br) \simeq H^{\text{Lie}}_*(\s_n ; \, \br) 
\ot [\Lambda^*(I_n)]^{\s_n}, \\
& H^{\text{Lie}}_*(\g_n ; \, \g_n) \simeq H^{\text{Lie}}_*(\s_n ; \, \br) 
\ot [\g_n \ot \Lambda^*(I_n)]^{\s_n}
\end{align*}
\end{lemma}
\begin{proof}
The lemma follows essentially from the Hochschild-Serre spectral
sequence \cite{HS}, the application of which we briefly outline to aid
in the identification of representative homology cycles, and to
reconcile the lemma with its cohomological version in \cite{HS}.  
Consider the filtration $\mathcal{F}_m$, $m \geq -1$, of the complex
$\Lambda^*(\g_n)$ given by
\begin{align*}
& {\mathcal{F}}_{-1} = \{ 0 \}, \\
& {\mathcal{F}}_0 = \Lambda^* (I_n), \ \ \ 
{\mathcal{F}}_0^k = I_n^{\w k},
  \ \ k = 0, \ 1, \ 2, \ 3, \ \ldots , \\
& {\mathcal{F}}_m^k = \{ g_1 \wedge \, \ldots \, \wedge g_{k+m} \in \g^{\wedge
    (k+m)}_n \ | \ {\text{at most $m$-many $g_i$'s}} \notin I_n \}.
\end{align*}
Then each $\mathcal{F}_m$ is a chain complex, and $\mathcal{F}_m$ is
a subcomplex of $\mathcal{F}_{m+1}$.  For $m \geq 0$, we have
$$ E_{m,\, k}^0 = {\mathcal{F}}_m^k/{\mathcal{F}}_{m-1}^k \simeq
I_n^{\wedge k} \ot (\g_n /I_n)^{\wedge m}.  $$
Since $I_n$ is Abelian and the action of $I_n$ on $\g_n /I_n$ is
trivial, it follows that
$$  E_{m, \, k}^1 \simeq I_n^{\wedge k} \ot (\g_n /I_n)^{\wedge m}. $$
Using the isomorphism $\g_n /I_n \simeq \s_n$, we have
$$  E_{m, \, k}^2 \simeq H_m(\s_n ; \, I_n^{\wedge k}). $$
Now, $\s_n$ is a simple Lie algebra and as an $\s_n$-module
$$  I_n^{\wedge k} \simeq (I_n^{\wedge k})^{\s_n} \oplus M, $$
where $M \simeq M_1 \oplus M_2 \oplus \ldots \oplus M_t$ is a
direct sum of simple modules on which $\s_n$ acts non-trivially.
Hence 
$$  H_*(\s_n ;\, I_n^{\wedge k}) \simeq H_*(\s_n ; \, (I_n^{\wedge
  k})^{\s_n}) \oplus H_*(\s_n ; \, M).  $$
Clearly,
\begin{align*}
& H_*(\s_n ; \, (I_n^{\wedge k})^{\s_n}) \simeq H_*(\s_n ; \, \br )
  \ot (I_n^{\wedge k})^{\s_n} \\
& H_*(\s_n ; \, M) \simeq \sum_{i=1}^t H_*(\s_n ; \, M_i) \simeq 0,
\end{align*}
where the latter isomorphism holds since each $M_i$ is simple with
non-trivial $\s_n$ action.  See \cite[Prop. VII.5.6]{Hilton} for more
details. 

Let $\theta$ be a cycle in $\Lambda^m (\s_n)$ representing an element
of $H_m(\s_n ; \, \br)$, and let $z \in (I_n^{\wedge
k})^{\s_n}$.  Then $z \w \theta \in \g_n^{\wedge (m+k)}$
represents an absolute cycle in $\Lambda^*(\g_n)$, since, if $\theta$
is a sum of elements of the form $s_1 \w s_2 \w \, \ldots \, \w s_m$,
then $[z, \ s_i] = 0$ for each $s_i \in \s_n$.     
Thus, $E_{m, \, k}^2 \simeq E_{m, \, k}^{\infty}$, and
$$  H_*(\g_n ; \, \br) \simeq 
H_*(\s_n ; \, \br) \ot [\Lambda^*(I_n)]^{\s_n}. $$ 
By a similar filtration and spectral sequence argument for 
$\g_n \ot \Lambda^*(\g_n)$, we have
$$  H_*(\g_n ; \, \g_n ) \simeq 
H_*(\s_n ; \, \br) \ot [\g_n \ot \Lambda^*(I_n)]^{\s_n}. $$  
\end{proof}

Let $\om_n = \sum_{i=1}^n \frac{\p}{\p x^i} \wedge \frac{\p}{\p y^i}
\in I_n^{\wedge 2}$.  One checks that $\om_n \in (I_n^{\wedge
2})^{\s_n}$ against the basis for $\s_n$ given in \S 2.  It follows
that
$$  \om_n^{\wedge k} \in [I_n^{\wedge 2k}]^{\s_n}.  $$
Letting $\Lambda^*(\om_n)$ denote the exterior algebra generated by
$\om_n$, we prove in the appendix that
\begin{lemma} \label{3.2}
There are isomorphisms
\begin{align*}
& [\Lambda^*(I_n)]^{\s_n} \simeq \Lambda^*(\om_n) := \sum_{k \geq 0}
  \Lambda^k (\om_n), \\
& [\g_n \ot \Lambda^*(I_n)]^{\s_n} \simeq \bar{\Lambda}^*(\om_n) 
:= \sum_{k \geq 1}  \Lambda^k (\om_n)
\end{align*}
where the first is an isomorphism of algebras, and the second is an
isomorphism of vector spaces.  
\end{lemma}

Combining this with Lemma \eqref{3.1}, we have
\begin{lemma} \label{3.3}
There are vector space isomorphisms
\begin{align*}
& H^{\text{Lie}}_* (\g_n ; \, \br) \simeq H_*(\s_n ; \, \br) \ot
  \Lambda^*(\om_n), \\
& H^{\text{Lie}}_* (\g_n ; \, \g_n) \simeq H_*(\s_n ; \, \br) \ot
  \bar{\Lambda}^*(\om_n).
\end{align*}
\end{lemma}

It is known that for cohomology,
$$ H^*_{\text{Lie}} (\s_n ; \, \br) \simeq \Lambda^*(u_3, \, u_7, \,
u_{11}, \, \ldots , \, u_{4n-1}), $$
where $u_i$ is a class in dimension $i$.  Also, 
$$  H^{\text{Lie}}_k (\s_n ; \, \br) \simeq 
H^k_{\text{Lie}} (\s_n ; \, \br).  $$
See the reference \cite[p.\ 343]{Whitehead} for the homology of the
symplectic Lie group.

\section{The Leibniz Homology of $\g_n$}

Recall that for a Lie algebra $\g$ over a ring $k$, and more generally
for a Leibniz algebra $\g$ \cite{LP}, the Leibniz homology of $\g$,
written $HL_*(\g ; \, k)$, is the homology of the chain complex
$T(\g)$:
$$  \CD 
k @<0<< \g @<[\ , \ ]<< \g^{\ot 2} @<<< \ldots @<<< \g^{\ot (n-1)}
@<d<< \g^{\ot n} @<<< \ldots,  
\endCD $$
where 
\begin{align*}
& d(g_1, \, g_2, \, \ldots , \,  g_n) = \\
& \sum_{1 \leq i < j \leq n} (-1)^j \, (g_1, \, g_2, \, \ldots, \, 
g_{i-1}, \, [g_i, \, g_j], \, g_{i+1}, \, \ldots \, 
\hat{g}_j \,\ldots , \, g_n), 
\end{align*}
and $(g_1, \, g_2, \, \ldots, \, g_n)$ denotes the element 
$g_1 \ot g_2 \ot \, \ldots \, \ot g_n \in \g^{\ot n}$.  

The canonical projection $\pi_1 : \g^{\ot n} \to \g^{\w n}$, $n \geq 0$, 
is a map of chain complexes, $T(\g ) \to \Lambda^*(\g )$,  
and induces a $k$-linear map on homology
$$  HL_*(\g ; \, k) \to H^{\text{Lie}}_*(\g ; \, k).  $$
Letting
$$  ({\text{ker}}\,\pi_1)_n [2] = {\text{ker}}\, [\g^{\ot (n+2)} \to
    \g^{\w (n+2)}], \ \ \ n \geq 0,  $$
Pirashvili \cite{Pirashvili} defines the relative theory
$H^{\text{rel}}(\g)$ as the homology of the complex
$$  C^{\text{rel}}_n (\g) = ({\text{ker}}\,\pi_1)_n [2],  $$
and studies the resulting long exact sequence relating Lie and Leibniz
homology:  
$$  \CD \label{4.1} 
\cdots @>{\p}>> H^{\text{rel}}_{n-2}(\g) @>>>   HL_n(\g) @>>>
H^{\text{Lie}}_n(\g) @>{\p}>>  H^{\text{rel}}_{n-3}(\g) @>>> \\
\cdots @>{\p}>> H^{\text{rel}}_0(\g) @>>>
HL_2(\g) @>>> H^{\text{Lie}}_2(\g) @>>> 0 \\
@.  0 @>>> HL_1(\g) @>>> H^{\text{Lie}}_1(\g) @>>> 0 \\
@.  0 @>>> HL_0(\g) @>>> H^{\text{Lie}}_0(\g) @>>> \phantom{.}0. 
\endCD $$

An additional exact sequence is required for calculations of $HL_*$.  
Consider the projection
$$  \pi_2 : \g \ot \g^{\w n} \to \g^{\w (n+1)}, \ \ \ n \geq 0,  $$
and the resulting chain map
$$  \pi_2 : \g \ot \Lambda^*(\g) \to \Lambda^{*+1}(\g).  $$
Let $HR_n(\g)$ denote the homology of the complex
$$  CR_n(\g) = ({\text{ker}}\, \pi_2)_n[1] = 
{\text{ker}}\, [\g \ot \g^{\w (n+1)} \to \g^{\w (n+2)}], \ \ \ n 
\geq 0.  $$
There is a resulting long exact sequence
$$  \CD
\cdots @>{\p}>> HR_{n-1}(\g) @>>> H^{\text{Lie}}_n(\g ; \, \g) @>>>
H^{\text{Lie}}_{n+1}(\g) @>{\p}>>  \\ 
\cdots @>{\p}>> HR_0(\g) @>>> H^{\text{Lie}}_1(\g ; \, \g) @>>>
H^{\text{Lie}}_2(\g) @>{\p}>>  \\ 
@.  0 @>>> H^{\text{Lie}}_0(\g ; \, \g) @>>>
H^{\text{Lie}}_1(\g) @>>> 0.
\endCD  $$
The projection $\pi_1 : \g^{\ot (n+1)} \to \g^{\w (n+1)}$ can be
written as the composition of projections
$$  \g^{\ot (n+1)} \longrightarrow \g \ot \g^{\w n}  
\longrightarrow \g^{\w (n+1)},  $$
which leads to a natural map between exact sequences
$$  \CD
H^{\text{rel}}_{n-1}(\g)  @>>> HL_{n+1}(\g) @>>>
H^{\text{Lie}}_{n+1}(\g)  @>{\p}>> H^{\text{rel}}_{n-2}(\g) \\
@VVV   @VVV   @V{\mathbf{1}}VV   @VVV  \\
HR_{n-1}(\g)  @>>> H^{\text{Lie}}_n(\g ; \, \g)  @>>>
H^{\text{Lie}}_{n+1}(\g)  @>{\p}>> HR_{n-2}(\g) 
\endCD  $$
and an articulation of their respective boundary maps $\p$.

\begin{lemma} \label{4.2}
For the affine symplectic Lie algebra $\g_n$, there is a natural
isomorphism
$$  H_k(\s_n ; \, \br) \overset{\simeq}{\longrightarrow}
HR_{k-3}(\g_n ; \, \br), \ \ \ k \geq 3, $$
that factors as the composition
$$ H_k^{\text{Lie}}( \s_n ;\, \br) 
\underset{\p}{\overset{\simeq}{\longrightarrow}}
HR_{k-3}(\s_n ; \, \br) \overset{\simeq}{\longrightarrow}
HR_{k-3}(\g_n ; \br), $$
and the latter isomorphism is induced by the inclusion $\s_n
\hookrightarrow \g_n$.  
\end{lemma}

\begin{proof}
Since $\s_n$ is a simple Lie algebra, from
\cite[Prop. VII.5.6]{Hilton} we have  
$$  H_k^{\text{Lie}}(\s_n ; \, \s_n) = 0 , \ \ \ k \geq 0.  $$
From the long exact sequence
$$ \ \cdots \longrightarrow HR_{k-1}(\s_n ; \, \br) \longrightarrow
H_k^{\text{Lie}}(\s_n; \, \s_n) \longrightarrow
H_{k+1}^{\text{Lie}}(\s_n; \, \br) \overset{\p}{\longrightarrow} 
\ \cdots \, ,$$
it follows that $\p : H_k^{\text{Lie}}(\s_n; \, \br) \to
HR_{k-3}(\s_n; \, \br)$ is an isomorphism for $k \geq 3$.  The 
inclusion of Lie algebras $\s_n \hookrightarrow \g_n$ induces a map of exact
sequences
$$ \CD
@>>> HR_{k-1}(\s_n; \, \br) @>>> H_k^{\text{Lie}}(\s_n; \, \s_n) @>>>
H_{k+1}^{\text{Lie}}(\s_n; \, \br) @>\p>>  \\
@. @VVV  @VVV  @VVV  \\
@>>> HR_{k-1}(\g_n; \, \br) @>>> H_k^{\text{Lie}}(\g_n; \, \g_n) @>>>
H_{k+1}^{\text{Lie}}(\g_n; \, \br) @>\p>> 
\endCD $$
From Lemma \eqref{3.3}
\begin{align*}
& H_*^{\text{Lie}}(\g_n ; \, \br) \simeq H_*(\s_n ; \, \br) \ot
  \Lambda^*(\om_n) \\
& H_*^{\text{Lie}}(\g_n ; \, \g_n) \simeq H_*(\s_n ; \, \br) \ot
  \bar{\Lambda}^*(\om_n) .
\end{align*}
The map 
$H_*^{\text{Lie}}(\g_n ; \, \g_n) \to H_{*+1}^{\text{Lie}}(\g_n ; \, \br)$
is an inclusion on homology with cokernel $H_{*+1}^{\text{Lie}}(\s_n ;
\, \br)$.  The result now follows from the map between exact sequences
and a knowledge of the generators of $H_*^{\text{Lie}}(\g_n ; \, \br)$
gleaned from Lemma \eqref{3.1}.  
\end{proof}
\begin{theorem} \label{4.3}
There is an isomorphism of vector spaces
$$  HL_*(\g_n ; \, \br) \simeq \Lambda^*(\om_n)  $$
and an algebra isomorphism
$$  HL^*(\g_n ; \, \br) \simeq \Lambda^*(\om_n^*), \ \ \ 
\om_n^* = \sum_{i=1}^n dx^i \w dy^i , $$
where $HL^*$ is afforded the shuffle algebra.
\end{theorem}
\begin{proof}
Consider the Pirashvili filtration \cite{Pirashvili} of the complex
$$  C_n^{\text{rel}}(\g) = {\text{ker}}( \g^{\ot(n+2)} \to
\g^{\w (n+2)}), \ \ \ n \geq 0,  $$
given by 
$$ \mathcal{F}_m^k(\g) = \g^{\ot k} \ot {\text{ker}}( \g^{\ot(m+2)} \to
\g^{\w (m+2)}), \ \ \ m \geq 0, \ k \geq 0. $$
Then $\mathcal{F}_m^*$ is a subcomplex of $\mathcal{F}_{m+1}^*$ and
the resulting spectral sequence converges to $H_*^{\text{rel}}(\g)$.
From \cite{Pirashvili} we have
$$  E^2_{m, \, k} \simeq HL_k(\g) \ot HR_m(\g), \ \ \ m \geq 0, \ k
\geq 0.  $$
From the proof of Lemma \eqref{4.2}, there is an isomorphism
$$  \p :  H_3^{\text{Lie}}( \g_n ; \, \br)
\overset{\simeq}{\longrightarrow} HR_0(\g_n ; \, \br) \simeq \br. $$
From the long exact sequence relating Lie and Leibniz homology, 
it follows that $HL_2(\g_n ; \, \br) \to H_2^{\text{Lie}}(\g_n ; \, \br)$ 
is an isomorphism.  Since
$$  \tilde{\om}_n = \frac{1}{2} \sum_{i=1}^n \bigg( \frac{\p}{\p x^i} \ot
\frac{\p}{\p y^i} - \frac{\p}{\p y^i} \ot \frac{\p}{\p x^i} \bigg)  $$
is a cycle in the Leibniz complex that maps to $\om_n$ in the Lie
algebra complex, it follows that $\tilde{\om}_n$ generates 
$HL_2(\g_n; \, \br)$.  

We claim that all elements in $HL_0(\g_n) \ot HR_*(\g_n)$ are absolute
cycles.  The inclusion $\s_n \hookrightarrow \g_n$ induces a map
between exact sequences
$$ \CD
HL_k(\s_n) @>>> H_k^{\text{Lie}}(\s_n) @>\p>>
H_{k-3}^{\text{rel}}(\s_n) @>>> HL_{k-1}(\s_n) \\
@VVV  @VVV  @VVV  @VVV  \\
HL_k(\g_n) @>>> H_k^{\text{Lie}}(\g_n) @>\p>>
H_{k-3}^{\text{rel}}(\g_n) @>>> HL_{k-1}(\g_n)
\endCD  $$
Since $\s_n$ is a simple Lie algebra, $HL_k(\s_n ; \, \br) = 0$, $k
\geq 1$ \cite{Ntolo}.  Thus, 
$\p : H_k^{\text{Lie}}(\s_n) \to H_{k-3}^{\text{rel}}(\s_n)$
is an isomorphism for $k \geq 3$.  The inclusion
$\mathcal{F}_m^*(\s_n) \hookrightarrow \mathcal{F}_m^*(\g_n)$ induces
a map of spectral sequences, and hence a map
$$  HL_0(\s_n) \ot HR_*(\s_n) \longrightarrow HL_0(\g_n) \ot
HR_*(\g_n).  $$
Since $HR_*(\s_n) \simeq H_*^{\text{rel}}(\s_n)$, all classes in
$HL_0(\s_n) \ot HR_*(\s_n)$ are absolute cycles.  Now, $HR_*(\s_n)$
maps isomorphically to $HR_*(\g_n)$, and by naturality, all classes in
$HL_0(\g_n) \ot HR_*(\g_n)$ are absolute cycles.  Moreover,
$$  \p : H_*^{\text{Lie}}(\g_n) \to H_{*-3}^{\text{rel}}(\g_n)  $$
maps the classes in $\bar{H}_*^{\text{Lie}}(\s_n)$ injectively to 
$H_{*-3}^{\text{rel}}(\g_n)$ in the diagram
$$ \CD
0 @>>> H_*^{\text{Lie}}(\s_n) @>\p>> H_{*-3}^{\text{rel}}(\s_n) @>>> 0 \\
@.  @VVV  @VVV  \\
\cdots @>>> H_*^{\text{Lie}}(\g_n) @>\p>> H_{*-3}^{\text{rel}}(\g_n)
@>>> \cdots 
\endCD  $$
where the vertical arrows are inclusions.

We claim that all elements in $HL_2(\g_n) \ot HR_*(\g_n)$ are absolute
cycles as well.  Let $[\theta] \in HR_m(\g_n)$ be represented by the
sum
$$  \theta = \sum_{j=1}^n X_{1, \, j} \ot X_{2, \, j} \w X_{3, \, j}
\w \ldots \w X_{m+1, \, j}\; ,  $$
where each $X_{i, \, j} \in \s_n$ and $d\theta = 0$.  By invariance,
$$  [\tilde{\om}_n, \, X_{i, \, j}] = 0 
\ \ \text{for each} \ \ X_{i, \, j}.  $$
It follows that $d(\tilde{\om}_n \ot \theta) = d(\tilde{\om}_n) \ot
\theta + \tilde{\om}_n \ot d\theta = 0$, and $\tilde{\om}_n \ot \theta$
represents an absolute cycle in $H_*^{\text{rel}}(\g_n)$.  To compute
$$  \p : H_*^{\text{Lie}}(\g_n) \to H_{*-3}^{\text{rel}}(\g_n)  $$
on classes of the form $[\om_n] \ot \bar{H}_*^{\text{Lie}}(\s_n)$, let
$[\theta'] \in \bar{H}_*^{\text{Lie}}(\s_n)$ with $\p (\theta') =
\theta$.  By lifting $\om_n \w \theta'$ to $\tilde{\om}_n \ot \theta'$
in $T(\g_n)$ and using invariance, we have
$$  \p (\om_n \w \theta') = \tilde{\om}_n \ot \p (\theta') = 
\tilde{\om}_n \ot \theta.  $$
At this point $H_k^{\text{rel}}(\g_n)$ is completely determined for $k
\leq 2$.  By an examination of $H_1^{\text{rel}}(\g_n)$, 
$$  \om_n^{\w 2} \in {\text{ker}} \, \p , \ \ \ 
\p : H_4^{\text{Lie}}(\g_n) \to H_1^{\text{rel}}(\g_n).  $$
Thus, $(\tilde{\om}_n)^{\w 2}$ generates a non-zero class in $HL_4(\g_n)$
mapping to the class $\om_n^{\w 2} \in H_4^{\text{Lie}}(\g_n)$.  As
before, all classes in $HL_4(\g_n) \ot HR_*(\g_n)$ are absolute cycles
and in Im$\, \p$.  Thus, $H_k^{\text{rel}}(\g_n)$ is completely
determined for $k \leq 4$.  By induction on $k$, $(\tilde{\om_n})^{\w k}$
is a non-zero class in $HL_{2k}(\g_n)$, and
\begin{align*}
& H_*^{\text{rel}}(\g_n) \simeq \Lambda^*(\om_n) \ot HR_*(\g_n) \simeq
  \Lambda^*(\om_n) \ot H_{*+3}^{\text{Lie}}(\s_n) \\
& HL_*(\g_n) \simeq \Lambda^*(\om_n) .
\end{align*}

For the cohomology isomorphism
$$  HL^*(\g_n ; \, \br) \simeq \Lambda^*(\om_n^*), \ \ \
\om_n^* = \sum_{i=1}^n dx^i \w dy^i, $$
where $dx^i$ is the dual of $\p x^i$ and $dy^i$ the dual of $\p y^i$
with respect to the basis of $\g_n$ given by $\mathcal{B}_1 \cup
\mathcal{B}_2$ given in \S 2.  Since
$$  HL^*(\g_n ; \, \br) \simeq {\text{Hom}}( HL_*(\g_n ; \, \br), \
\br), $$
the result follows by using the full shuffle product on cochains.
\end{proof}

\section{A Characteristic Map}

Let $M$ be a symplectic manifold, $\x (M)$ the Lie algebra of $C^{\infty}$ 
vector fields on $M$, and $\x_{H}(M)$ the Lie algebra of Hamiltonian vector 
fields \cite[p. 85]{MS}.  From \cite{Lodder}, there is a natural map
$$  H^*_{dR}(M; \, \br) \longrightarrow HL^*(\x (M); \, C^{\infty}(M)),  $$
where $H^*_{dR}(M)$ denotes deRham cohomology.  The inclusion of Lie algebras
$\x_{H}(M) \hookrightarrow \x (M)$ induces a (contravariant) map
$$  HL^*(\x (M); \, C^{\infty}(M)) \longrightarrow
HL^*(\x_{H}(M); \, C^{\infty}(M))  $$
on cohomology, while the inclusion of coefficients $\br \to C^{\infty}(M)$ 
induces a (covariant) map
$$  HL^*(\x_{H}(M); \, \br) \longrightarrow HL^*(\x_{H}(M); \, C^{\infty}(M)). $$  
Let $p \in M$ and let $U$ be an open neighborhood of $p$ homeomorphic to 
$R^{2n}$ in the atlas of charts for $M$.  There is a natural morphism of Lie 
algebras $\x_{H}(M) \to \x_{H}(U)$ given by the restriction of vector fields from
$M$ to $U$, and resulting linear maps
$$ HL^*(\x_{H}(U); \, \br) \to HL^*(\x_{H}(M); \br)  \to 
HL^*(\x_H (M);\, C^{\infty}(M)).  $$ 
Now, $\x_{H}(U) \simeq \x_{H}(\br^{2n})$ as Lie algebras, and thus 
there are local  maps
$$  HL^*(\x_{H}(\br^{2n}); \, \br)(p) \to HL^*(\x_{H}(M); \, \br) $$
for each $p \in M$.  Finally, the inclusion
$\g_n \hookrightarrow \x_{H}(\br^{2n})$ is a morphism of Lie algebras
inducing a map on Leibniz homology
$$  HL_*(\g_n; \, \br) \longrightarrow HL_*(\x_{H}(\br^{2n}); \, \br). $$

\section{Appendix}

The goal of the appendix is to establish Lemma \eqref{3.2}, namely the vector
space isomorphisms
\begin{align}
& [\Lambda^*(I_n)]^{\s_n} \simeq \Lambda^*(\om_n) := \sum_{k \geq 0}
  \Lambda^k (\om_n) \label{6.1} \\
& [\g_n \ot \Lambda^*(I_n)]^{\s_n} \simeq \bar{\Lambda}^*(\om_n)
:= \sum_{k \geq 1}  \Lambda^k (\om_n), \label{6.2}
\end{align}
where the former is also an algebra isomorphism.  First, note that as
an $\s_n$-module, $\g_n \simeq I_n \oplus \s_n$, and
$$  [\g_n \ot \Lambda^*(I_n)]^{\s_n} \simeq 
[I_n \ot \Lambda^*(I_n)]^{\s_n} \oplus [\s_n \ot \Lambda^*(I_n)]^{\s_n}. $$
Thus, line \eqref{6.2} would follow from the vector space isomorphisms
\begin{align*}
& [I_n \ot \Lambda^*(I_n)]^{\s_n} \simeq \bar{\Lambda}^*(\om_n)  \\
& [\s_n \ot \Lambda^*(I_n)]^{\s_n} = \{ 0 \}.  
\end{align*}

We first demonstrate isomorphism \eqref{6.1} in the following lemma.
\begin{lemma}
$$  [\Lambda^*(I_n)]^{\s_n} \simeq \Lambda^*(\om_n).  $$ 
\end{lemma}
\begin{proof}
We proceed by induction on $n$.  For $n = 1$, 
\begin{align*}
& I_1 = \bigg{\langle} \frac{\p}{\p x^1}, \ \frac{\p}{\p y^1} \bigg{\rangle} \\ 
& \s_1 = \bigg{\langle} x_1 \frac{\p}{\p y^1}, \ y_1 \frac{\p}{\p x^1}, \ 
y_1 \frac{\p}{\p y^1} - x_1 \frac{\p}{\p x^1} \bigg{\rangle}. 
\end{align*}
By direct calculation, $(I_1)^{\s_1} = \{ 0 \}$, and 
$(I_1^{\w 2})^{\s_1} = \langle \frac{\p}{\p x^1} \w \frac{\p}{\p y^1} \rangle$.  

By the inductive hypothesis, suppose
$$  [\Lambda^*(I_{n-1})]^{\s_{n-1}} = \Lambda^*(\om_{n-1}).  $$
Consider then two cases for $I_n^{\w k}$, $k$ odd, and $k$ even.  For $k$ odd,
let $z \in I_n^{\w k}$ and consider
$$  z = z_1 + z_2 \w \frac{\p}{\p x^n} + z_3 \w \frac{\p}{\p y^n} +
z_4 \w \frac{\p}{\p x^n} \w \frac{\p}{\p y^n},  $$
where $z_1 \in I_{n-1}^{\w k}$, $z_2$, $z_3 \in I_{n-1}^{\w (k-1)}$,
and $z_4 \in I_{n-1}^{\w (k-2)}$.  Note that 
$$  [z, \ y_n \frac{\p}{\p y^n} - x_n \frac{\p}{\p x^n}] = 
- z_2 \w \frac{\p}{\p x^n} + z_3 \frac{\p}{\p y^n}.  $$
For $z \in (I_n^{\w k})^{\s_n}$, 
$[z, \ y_n \frac{\p}{\p y^n} - x_n \frac{\p}{\p x^n}] = 0$, and
$$  z = z_1 + z_4 \w \frac{\p}{\p x^n} \w \frac{\p}{\p x^n}.  $$
For any $X \in \s_{n-1} \subseteq \s_n$, we have 
$$  0 = [z, \ X] = [z_1, \ X] + 
[z_4, \ X] \w \frac{\p}{\p x^n} \w \frac{\p}{\p y^n}.  $$
If non-zero, the terms $[z_1, \ X]$ and 
$[z_4, \ X] \w \frac{\p}{\p x^n} \w \frac{\p}{\p y^n}$ are linearly
independent and would not sum to zero.  Thus,
$$  z_1 \in (I_{n-1}^{\w k})^{\s_{n-1}} = \{ 0 \}, \ \ \ 
z_4 \in (I_{n-1}^{\w (k-2)})^{\s_{n-1}} = \{ 0 \}.  $$
It follows that $(I_n^{\w k})^{\s_n} = \{ 0 \}$ for $k$ odd.   

For $k$ even, let $k = 2q$, $z \in (I_n^{\w 2q})^{\s_n}$, and repeat the above
argument to the point
\begin{align*}
& z_1 \in (I_{n-1}^{\w 2q})^{\s_{n-1}} = \langle \om_{n-1}^{\w q} \rangle \\
& z_4 \in (I_{n-1}^{\w 2(q-1)})^{\s_{n-1}} = 
\langle \om_{n-1}^{\w (q-1)} \rangle 
\end{align*}
Thus, $z = c_1 \om_{n-1}^{\w q} + c_2 \om_{n-1}^{\w (q-1)} \w
\frac{\p}{\p x^n} \w \frac{\p}{\p y^n}$, $c_1$, $c_2 \in \br$.  Bracketing
with $X = x_1 \frac{\p}{\p y^n} + x_n \frac{\p}{\p y^1}$ yields 
$$  0 = [z, \ X] = (c_2 - qc_1)\om_{n-1}^{\w (q-1)} \w 
\frac{\p}{\p y^1} \w \frac{\p}{\p y^n}.  $$
Hence, $z$ is a real multiple of
$$  \om_{n-1}^{\w q} + q \om_{n-1}^{\w (q-1)} 
\w \frac{\p}{\p x^n} \w \frac{\p}{\p y^n} =  
{\bigg( \om_{n-1} + \frac{\p}{\p x^n} \w \frac{\p}{\p y^n}\bigg)}^{\w q} =  
\om_n^{\w q}.  $$
\end{proof}

\begin{lemma}
$$ [I_n \ot \Lambda^*(I_n)]^{\s_n} \simeq \bar{\Lambda}^*(\om_n). $$
\end{lemma}
\begin{proof}
The proof proceeds by induction on $n$.  For $n = 1$, a direct verification
yields
$$  (I_1)^{\s_1} = \{ 0 \}, \ \ \ 
(I_1 \ot I_1)^{\s_1} = \bigg{\langle} \frac{\p}{\p x^1} \w \frac{\p}{\p y^1}
\bigg{\rangle},  $$
where $\frac{\p}{\p x^1} \w \frac{\p}{\p y^1} = \frac{\p}{\p x^1} \ot \frac{\p}{\p y^1}
- \frac{\p}{\p y^1} \ot \frac{\p}{\p x^1}.$  Also, 
$(I_1 \ot I_1^{\w 2})^{\s_1} = \{ 0 \}$ by direct calculation.   
The inductive hypothesis states 
$$ [I_{n-1} \ot \Lambda^*(I_{n-1})]^{\s_{n-1}} \simeq 
\bar{\Lambda}^*(\om_{n-1}).  $$
Let $v \in I_n \ot I_n^{\w k}$, $v = u_1 + u_2$, where 
$$  u_1 \in I_{n-1} \ot I_{n-1}^{\w k}, \ \ \ 
u_2 \in (I_n \ot I_n^{\w k})/(I_{n-1} \ot I_{n-1}^{\w k}).  $$
A vector space basis of $(I_n \ot I_n^{\w k})/(I_{n-1} \ot I_{n-1}^k)$ is
given by the families of elements:
\begin{enumerate}
\item[(1)] $\frac{\p}{\p x^n} \ot \frac{\p}{\p x^n} \w \frac{\p}{\p y^n} \w 
\frac{\p}{\p z^1} \w \frac{\p}{\p z^2} \w \ldots \w \frac{\p}{\p z^{k-2}}$

\item[(2)]   $\frac{\p}{\p y^n} \ot \frac{\p}{\p x^n} \w \frac{\p}{\p y^n} \w 
\frac{\p}{\p z^1} \w \frac{\p}{\p z^2} \w \ldots \w \frac{\p}{\p z^{k-2}}$

\item[(3)] $\frac{\p}{\p x^n} \ot \frac{\p}{\p x^n} \w 
\frac{\p}{\p z^1} \w \frac{\p}{\p z^2} \w \ldots \w \frac{\p}{\p z^{k-1}}$

\item[(4)] $\frac{\p}{\p x^n} \ot \frac{\p}{\p y^n} \w 
\frac{\p}{\p z^1} \w \frac{\p}{\p z^2} \w \ldots \w \frac{\p}{\p z^{k-1}}$

\item[(5)] $\frac{\p}{\p y^n} \ot \frac{\p}{\p x^n} \w 
\frac{\p}{\p z^1} \w \frac{\p}{\p z^2} \w \ldots \w \frac{\p}{\p z^{k-1}}$

\item[(6)] $\frac{\p}{\p y^n} \ot \frac{\p}{\p y^n} \w 
\frac{\p}{\p z^1} \w \frac{\p}{\p z^2} \w \ldots \w \frac{\p}{\p z^{k-1}}$

\item[(7)] $\frac{\p}{\p z^1} \ot \frac{\p}{\p x^n} \w 
\frac{\p}{\p y^n} \w \frac{\p}{\p z^2} \w \ldots \w \frac{\p}{\p z^{k-1}}$

\item[(8)] $\frac{\p}{\p z^1} \ot \frac{\p}{\p x^n} \w 
\frac{\p}{\p z^2} \w \frac{\p}{\p z^3} \w \ldots \w \frac{\p}{\p z^{k}}$

\item[(9)] $\frac{\p}{\p z^1} \ot \frac{\p}{\p y^n} \w 
\frac{\p}{\p z^2} \w \frac{\p}{\p z^3} \w \ldots \w \frac{\p}{\p z^{k}}$,
\end{enumerate}
where, for each family, the $z^i$'s are elements of
$$  \{ x^1, \ x^2, \ \ldots , \ x^{n-1}, \ y^1, \ y^2, \ \ldots , \ y^{n-1} \}. $$   
Let $v \in (I_n \ot I_n^{\w k})^{\s_n}$ and $X = y_n \frac{\p}{\p y^n} - 
x_n \frac{\p}{\p x^n}$.  Then
$$  0 = [v, \ X] = [u_1 + u_2, \ X] = [u_2, \ X].  $$
To compute the $\s_n$-invariants, consider $u_2 \in {\text{ker(ad}}_X)$, where
${\text{ad}}_X (w) = [w, \ X]$.  The families (4), (5) and (7) above fall into
${\text{ker(ad}}_X)$.  Now consider $X = x_n \frac{\p}{\p y^n}$.  Family (7) along
with 
$$ \frac{\p}{\p x^n} \ot \frac{\p}{\p y^n} \w \frac{\p}{\p z^1} \w 
\ldots \w \frac{\p}{\p z^{k-1}}
- \frac{\p}{\p y^n} \ot \frac{\p}{\p x^n} \w \frac{\p}{\p z^1} \w 
\ldots \w \frac{\p}{\p z^{k-1}} $$
are elements of ${\text{ker(ad}}_X)$, $X = x_n \frac{\p}{\p y^n}$.  Then
$v = u_1 + s_1 + s_2$,
\begin{align*}
& s_1 = \sum_{z^1, \ldots, z^{k-1}} c_{1, \, *} \bigg( 
\frac{\p}{\p x^n} \ot \frac{\p}{\p y^n} \w \frac{\p}{\p z^1} \w 
\ldots \w \frac{\p}{\p z^{k-1}} \\
& \phantom{ s_1 = \sum_{z^1, \ldots, z^{k-1}} c_{1,*} \bigg(} 
- \frac{\p}{\p y^n} \ot \frac{\p}{\p x^n} \w \frac{\p}{\p z^1} \w 
\ldots \w \frac{\p}{\p z^{k-1}} \bigg) \\
& s_2 = \sum_{z^1, \ldots, z^{k-1}} c_{2,\, *} \bigg(
\frac{\p}{\p z^1} \ot \frac{\p}{\p x^n} \w 
\frac{\p}{\p y^n} \w \frac{\p}{\p z^2} \w \ldots \w \frac{\p}{\p z^{k-1}}
\bigg)
\end{align*} 
For $X \in \s_{n-1}$,
$$  0 = [v, \ X] = [u_1, \ X] + [s_1, \ X] + [s_2, \ X].  $$
Note that 
$$ [u_1, \ X] \in I_{n-1} \ot I_{n-1}^{\w k}, \ \ \
   [s_1, \ X] \notin I_{n-1} \ot I_{n-1}^{\w k}, \ \ \
   [s_2, \ X] \notin I_{n-1} \ot I_{n-1}^{\w k}.  $$
If non-zero, the summands of $[s_1, \ X]$ and $[s_2, \ X]$ would be
linearly independent.  Thus, $[s_1, \ X] = 0$, $[s_2, \ X] = 0$, and
$u_1 \in (I_{n-1} \ot I_{n-1}^{\w k})^{\s_{n-1}}$.  For $k$ even,
$(I_{n-1} \ot I_{n-1}^{\w k})^{\s_{n-1}} = \{ 0 \}$, $u_1 = 0$,
\begin{align*}
& [s_2, \ X] = \sum_{z^1, \ldots, z^{k-1}} c_{2, \, *}
[ \frac{\p}{\p z^1} \ot \frac{\p}{\p z^2} \w \ldots \w  
\frac{\p}{\p z^{k-1}}, \ X] \w \frac{\p}{\p x^n} \w \frac{\p}{\p y^n}, \\
& \sum_{z^1, \ldots, z^{k-1}} c_{2, \, *}
\frac{\p}{\p z^1} \ot \frac{\p}{\p z^2} \w \ldots \w  
\frac{\p}{\p z^{k-1}} \in (I_{n-1} \ot I_{n-1}^{\w (k-2)})^{\s_{n-1}}
= \{ 0 \}.
\end{align*}
Thus, $v = s_1$.  From 
$$ 0 = [s_1, \ x_n \frac{\p}{\p y^i} + x_i \frac{\p}{\p y^n}], \ \ \ 
0 = [s_1, \ y_i \frac{\p}{\p x^n} + y_n \frac{\p}{\p x^i}], $$
for $1 \leq i \leq n-1$, it follows that $s_1 = 0$.  

For $k$ odd, let $k = 2q - 1$.  Then
\begin{align*}
& u_1 \in (I_{n-1} \ot I_{n-1}^{\w k})^{\s_{n-1}} = 
\langle \om_{n-1}^{\w q} \rangle  \\
& \theta := \sum_{z^1, \ldots, z^{k-1}} c_{2, \, *}
\frac{\p}{\p z^1} \ot \frac{\p}{\p z^2} \w \ldots \w  
\frac{\p}{\p z^{k-1}} \in (I_{n-1} \ot I_{n-1}^{\w (k-2)})^{\s_{n-1}} 
= \langle \om_{n-1}^{\w (q-1)} \rangle  \\  
& u_1 = \lambda_1 \om_{n-1}^{\w q}, \ \ \ 
\theta = \lambda_2 \om_{n-1}^{\w (q-1)}, \ \ \ \lambda_1, \
\lambda_2 \in \br .
\end{align*} 
Note that 
\begin{align*}
& [ \lambda_1 \om_{n-1}^{\w q} + \lambda_2 \om_{n-1}^{\w (q-1)} \w
\frac{\p}{\p x^n} \w \frac{\p}{\p y^n}, \  y_i \frac{\p}{\p x^n} + y_n
\frac{\p}{\p x^i} ]  \\
& = (q \lambda_1 - \lambda_2) \om_n^{\w (q-1)} \w 
\frac{\p}{\p x^i} \w \frac{\p}{\p x^n}.  
\end{align*}
From
$$  0 = [ \lambda_1 \om_{n-1}^{\w q} + \lambda_2 \om_{n-1}^{\w (q-1)} \w
\frac{\p}{\p x^n} \w \frac{\p}{\p y^n} + s_1, \ X ]  $$
for $X = x_n \frac{\p}{\p y^i} + x_i \frac{\p}{\p y^n}$, 
$X = y_i \frac{\p}{\p x^n} + y_n \frac{\p}{\p x^i}$, $1 \leq i \leq
n-1$, it follows that $s_1 = 0$, and $(q \lambda_1 - \lambda_2) = 0$.
Letting $\lambda_1 = 1$, we have $\lambda_2 = q$, and 
$$  v = \om_{n-1}^{\w q} + q \om_{n-1}^{\w q} \w 
\frac{\p}{\p x^n} \w \frac{\p}{\p y^n} = \om_n^{\w q}.  $$
\end{proof}

\begin{lemma}
$$  [\s_n \ot \Lambda^*(I_n)]^{\s_n} = \{ 0 \}.  $$
\end{lemma}
\begin{proof}
We apply induction on $n$.  For $n = 1$, write a general element of
$\s_1 \ot \Lambda^*(I_1)$ as a linear combination of the basis elements
given in $\mathcal{B}_1$ and $\mathcal{B}_2$ of \S 2 ($n = 1$).  Then
apply $\text{ad}_X$ for $X = (y_1 \frac{\p}{\p y^1} - 
x_1 \frac{\p}{\p x^1})$.  The result $[\s_1 \ot \Lambda^*(I_1)]^{\s_1}
= \{ 0 \}$ follows from linear algebra.

Suppose that $[\s_{n-1} \ot \Lambda^*(I_{n-1})]^{\s_{n-1}} = \{ 0
\}$.  Since $\s_n$ is a simple Lie algebra, we have $(\s_n)^{\s_n} =
\{ 0 \}$.  Let $\mathcal{B}_1$ be the vector space basis for
$\s_{n-1}$ given in \S 2, and let 
\begin{align*}
& S = \{ x^1,\ x^2,\ \ldots, \ x^n, \ y^1, \ y^2, \ \ldots ,\ y^n \} \\ 
& S'= \{ x^1,\ x^2,\ \ldots,\ x^{n-1},\ y^1,\ y^2,\ \ldots ,\ y^{n-1} \}.
\end{align*}
A vector space basis of $(\s_n \ot I_n^{\w k})/(\s_{n-1} \ot
I_{n-1}^{\w k})$ is given by the families of elements:
\begin{enumerate}
\item[(1)] $e \ot \frac{\p}{\p x^n} \w \frac{\p}{\p y^n} \w 
\frac{\p}{\p z^1} \w \frac{\p}{\p z^2} \w \ldots \w \frac{\p}{\p z^{k-2}},
\ \ \ e \in \mathcal{B}_1, \ z^i \in S'$

\item[(2)] $e \ot \frac{\p}{\p x^n} \w \frac{\p}{\p z^1} 
\w \frac{\p}{\p z^2} \w \ldots \w \frac{\p}{\p z^{k-1}},
\ \ \ e \in \mathcal{B}_1, \ z^i \in S'$

\item[(3)] $e \ot \frac{\p}{\p y^n} \w \frac{\p}{\p z^1} 
\w \frac{\p}{\p z^2} \w \ldots \w \frac{\p}{\p z^{k-1}},
\ \ \ e \in \mathcal{B}_1, \ z^i \in S'$

\item[(4)] $(x_n  \frac{\p}{\p y^n}) \ot \frac{\p}{\p z^1} 
\w \frac{\p}{\p z^2} \w \ldots \w \frac{\p}{\p z^{k}},
\ \ \  z^i \in S$

\item[(5)] $( x_n  \frac{\p}{\p y^i} + x_i \frac{\p}{\p y^n}) \ot
\frac{\p}{\p z^1} \w \frac{\p}{\p z^2} \w \ldots \w \frac{\p}{\p z^{k}},
\ \ \  i < n, \ z^j \in S$

\item[(6)] $( y_n \frac{\p}{\p x^n}) \ot 
\frac{\p}{\p z^1} \w \frac{\p}{\p z^2} \w \ldots \w \frac{\p}{\p z^{k}},
\ \ \   z^i \in S$

\item[(7)] $( y_i  \frac{\p}{\p x^n} + y_n \frac{\p}{\p x^i}) \ot 
\frac{\p}{\p z^1} \w \frac{\p}{\p z^2} \w \ldots \w \frac{\p}{\p z^{k}},
\ \ \  i < n, \ z^j \in S$

\item[(8)] $( y_n \frac{\p}{\p y^n} - x_n \frac{\p}{\p x^n}) \ot 
\frac{\p}{\p z^1} \w \frac{\p}{\p z^2} \w \ldots \w \frac{\p}{\p z^{k}},
\ \ \   z^i \in S$

\item[(9)] $( y_i  \frac{\p}{\p y^n} - x_n \frac{\p}{\p x^i}) \ot 
\frac{\p}{\p z^1} \w \frac{\p}{\p z^2} \w \ldots \w \frac{\p}{\p z^{k}},
\ \ \  i < n, \ z^j \in S$

\item[(10)] $( y_n  \frac{\p}{\p y^i} - x_i \frac{\p}{\p x^n}) \ot 
\frac{\p}{\p z^1} \w \frac{\p}{\p z^2} \w \ldots \w \frac{\p}{\p z^{k}},
\ \ \  i < n, \ z^j \in S$
\end{enumerate}

Given $w \in (\s_n \ot I_n^{\w k})^{\s_n}$, let $w = u + v$, where 
$$ u \in (\s_{n-1} \ot I_{n-1}^{\w k}), \ \ \ 
v \in (\s_n \ot I_n^{\w k})/(\s_{n-1} \ot I_{n-1}^{\w k}).  $$
For all $X \in \s_n$,
$0 = \text{ad}_X (w) = \text{ad}_X (u) + \text{ad}_X (v)$.
Restricting to $X \in \s_{n-1}$, notice that if non-zero, the elements
$\text{ad}_X (u)$ and $\text{ad}_X (v)$ are linearly independent.
Thus, $\text{ad}_X (u) = 0$, and 
$$  u \in (\s_{n-1} \ot I_{n-1}^{\w k})^{\s_{n-1}} = \{ 0 \}.  $$
Now, $v$ can be written as a linear combination of the elements in
families (1)--(10).  We prove that $v = 0$ by applying the condition
$\text{ad}_X (v) = 0$ for successive choices of $X \in \s_n$.  First
apply $X = (y_n \frac{\p}{\p y^n} - x_n \frac{\p}{\p x^n})$, then $X
\in \s_{n-1}$ together with the inductive hypothesis.  Third, apply $X
= x_n \frac{\p}{\p y^n}$, fourth 
$X = (y_i \frac{\p}{\p y^i} - x_i \frac{\p}{\p x^i})$, fifth $X = x_i
\frac{\p}{\p y^i}$, and finally 
$X = (x_n \frac{\p}{\p y^i} + x_i \frac{\p}{\p y^n})$, where 
$1 \leq i \leq n-1$.  
\end{proof}

\end{document}